\documentclass[12pt,a4paper]{article}
\usepackage[utf8]{inputenc}
\usepackage{amsmath, amssymb, amsthm}
\usepackage{mathrsfs}
\usepackage{enumitem}

\theoremstyle{definition}
\newtheorem{definition}{Definition}[section]

\theoremstyle{plain}
\newtheorem{theorem}[definition]{Theorem}
\newtheorem{proposition}[definition]{Proposition}
\newtheorem{example}[definition]{Example}
\newtheorem{remark}[definition]{Remark}

\title{On Slender Generalized Groups}
\author{%
Mohammad Reza Ahmadi Zand\textsuperscript{1} \and 
Hamid Torabi Ardakani\textsuperscript{2}
}
\date{}

\begin{document}

\maketitle

\begin{center}
\footnotesize
\textsuperscript{1}Department of Mathematics, Yazd University, Yazd, Iran \\
Email: m.ahmadi@yazd.ac.ir \\
\textsuperscript{2}Department of Mathematics, Ferdowsi University of Mashhad, Mashhad, Iran \\
Email: hamid.torabi@um.ac.ir
\end{center}

\begin{abstract}
This paper introduces the concept of \emph{slender generalized groups}, extending the classical notion of slender abelian groups to the setting of generalized groups (also known as completely simple semigroups). We establish fundamental properties of slender generalized groups and prove that, in the abelian case, the classical and generalized definitions of slenderness coincide. Several structural results are provided, including the behavior of slenderness under homomorphisms and subgroups. We also present examples and non-examples to illustrate the theory. Our results show that slenderness is preserved under taking generalized subgroups and that the group components \( G_{e(a)} \) of a slender generalized group are slender in the classical sense.
\end{abstract}


\section{Introduction}
An interesting generalization of groups are generalized groups or completely simple semigroups \cite{Molaei1999GeneralizedGroups} which was introduced by M. Molaei for the first time.
A generalized group \cite{Molaei1999GeneralizedGroups} is a semigroup \( G \) which satisfies the following conditions:

(1) For each \( x \in G \) there exists a unique element \( z \in G \) such that \( zx = xz = x \) (we denote \( z \) by \( e(x) \)).

(2) For each \( x \in G \) there exists an element \( y \in G \) such that \( xy = yx = e(x) \).

It is well known that for a given \( x \in G \), \( e(e(x)) = e(x) \), and each \( x \) in \( G \) has a unique inverse in \( G \). The inverse of \( x \) is denoted by \( x^{-1} \) and \( (x^{-1})^{-1} = x \) \cite{Molaei1999GeneralizedGroups}. Moreover, \( e(x^{-1}) = e(x) \). Adeniran et al \cite{Adeniran2009}, show that an abelian generalized group is a group.

\begin{definition}[\cite{Molaei1999GeneralizedGroups}]
If \( G \) and \( H \) are two generalized groups, then a map \( f : G \rightarrow H \) is called a homomorphism if \( f(ab) = f(a)f(b) \) for all \( a, b \in G \).
\end{definition}

\begin{theorem}[\cite{Molaei1999GeneralizedGroups}]
Let \( f : G \rightarrow H \) be a homomorphism where \( G \) and \( H \) are two generalized groups. Then:

(1) \( f(e(a)) = e(f(a)) \),

(2) \( f(a^{-1}) = (f(a))^{-1} \),

for all \( a \in G \).
\end{theorem}

\begin{definition}[\cite{Molaei1999GeneralizedGroups}, \cite{Molaei2000Topological}]
If \( G \) is a generalized group and \( e(xy) = e(x)e(y) \) for all \( x, y \in G \), then \( G \) is called a normal generalized group.
\end{definition}

Let \( n \in \mathbb{N} \). The following notation is needed in what follows. The element \( i_n \in \prod_{\mathbb{N}} \mathbb{Z} \) is defined by \( i_n(k) = 1 \) if \( k = n \) and \( i_n(k) = 0 \) for all \( k \neq n \). We recall that an abelian group \( A \) is slender if for every homomorphism \( h : \prod_{\mathbb{N}} \mathbb{Z} \rightarrow A \), then \(\{n \in \mathbb{N} \mid h(i_n) \neq 0\}\) is a finite set. An abelian group \( A \) is slender if and only if for every homomorphism \( h : \prod_{\mathbb{N}} \mathbb{Z} \rightarrow A \), \(h(\prod_{\mathbb{N}\setminus\{1,2,\ldots,n\}}\mathbb{Z})=0\) for some \(n\in\mathbb{N}\). For example, every free abelian group is slender, however, divisible groups (e.g \(\mathbb{Q}\), \(J_{p}\) the p-adic integer group) are not slender. The subgroups and the direct sum of slender groups are slender \cite{Fuchs1970}. However, the direct product of slender groups is not slender in general. For example, the free abelian group \(\mathbb{Z}\) is slender and \(\prod_{\mathbb{N}}\mathbb{Z}\) is not slender, since the identity map \(1_{\prod_{\mathbb{N}}\mathbb{Z}}:\prod_{\mathbb{N}}\mathbb{Z}\to\prod_{\mathbb{N }}\mathbb{Z}\) is a homomorphism, which does not send any \(i_{n}\) to \(0\).

\begin{example}
The set \(\prod_{\mathbb{N}}\mathbb{Z}\) is a generalized group with the following operation.
\[
(x_{1},x_{2},x_{3},x_{4},x_{5},x_{6},\ldots)*(y_{1},y_{2},y_{3},y_{4},y_{5},y_{6}, \ldots)=(x_{1},y_{2},x_{3}+y_{3},x_{4},y_{5},x_{6}+y_{6},\ldots)
\]
Note that \(e(x_{1},x_{2},x_{3},x_{4},x_{5},x_{6}\ldots)=(x_{1},x_{2},0,x_{4},x_{5},0\ldots)\). We denote this generalized group by \(\prod_{\mathbb{N}}^{g}\mathbb{Z}\).
\end{example}

\begin{remark}
Recall that a map \(f:G_{1}\to G_{2}\) of generalized groups is called a homomorphism of generalized groups if \(f(ab)=f(a)f(b)\) for every \(a,b\in G_{1}\). For example, a map \(f:\prod_{\mathbb{N}}\mathbb{Z}\to\prod_{\mathbb{N}}^{g}\mathbb{Z}\) with \(f(x_{1},x_{2},x_{3},x_{4},x_{5},x_{6},\ldots)=(0,0,x_{3},0,0,x_{6},\ldots)\) is a homomorphism. Also a map \(g:\prod_{\mathbb{N}}^{g}\mathbb{Z}\to\prod_{\mathbb{N}}\mathbb{Z}\) with \(g(x_{1},x_{2},x_{3},x_{4},x_{5},x_{6},\ldots)=(x_{3},x_{6},x_{9},\ldots)\) is a homomorphism. But the identity map \(I:\prod_{\mathbb{N}}^{g}\mathbb{Z}\to\prod_{\mathbb{N}}\mathbb{Z}\) with \(I(x)=x\) is not a homomorphism.
\end{remark}

\begin{definition}
We call a generalized group \(G\) slender if for every homomorphism \(h:\prod_{\mathbb{N}}^{g}\mathbb{Z}\to G\) then \(\{n\in\mathbb{N} \mid h(i_n)\neq e(h(i_n))\}\) is a finite set.
\end{definition}

The proof of the following results is straightforward.

\begin{proposition}
If \(G\) is isomorphic to a generalized subgroup of a slender generalized group, then \(G\) is a slender generalized group.
\end{proposition}

By the above proposition, every generalized subgroup of a slender generalized group is a slender generalized group. Let \(G\) be a generalized slender group and \(a\in G\). Then, \(G_{e(a)}=\{g\in G \mid e(g)=e(a)\}\) with the operation of \(G\) is a group \cite{Molaei1999GeneralizedGroups}. Thus by Proposition 1.7, we have the following result.

\begin{proposition}
If \(G\) is a slender generalized group and \(a\in G\), then the group \(G_{e(a)}\) is a slender group.
\end{proposition}

\begin{example}
The generalized group \(\prod_{\mathbb{N}}^{g}\mathbb{Z}\) is not slender, since the identity map \(1_{\prod_{\mathbb{N}}^{g}\mathbb{Z}}:\prod_{\mathbb{N}}^{g}\mathbb{Z}\to\prod_{\mathbb{N}}^{g}\mathbb{Z}\) is a homomorphism, which does not send any \(i_{3n}\) to \(e(i_{3n})=0\).
\end{example}

In the following we show that an abelian group \(A\) is a slender group if and only if it is a slender generalized group.

\begin{theorem}
Let \(A\) be an abelian group. If \(A\) is a slender group then it is a slender generalized group.
\end{theorem}

\begin{proof}
Let \(h:\prod_{{\mathbb{N}}}^{g}{\mathbb{Z}}\to A\) be a homomorphism of generalized groups. Let \(f:\prod_{{\mathbb{N}}}{\mathbb{Z}}\to\prod_{{\mathbb{N}}}^{g}{\mathbb{Z}}\) be a homomorphism with \(f(x_{1},x_{2},x_{3},x_{4},x_{5},x_{6},\dots)=(0,0,x_{3},0,0,x_{6},\dots)\). Then \(h \circ f:\prod_{{\mathbb{N}}}{\mathbb{Z}}\to A\) is a homomorphism of abelian groups. Since \(A\) is a slender group, \(\{n\in{\mathbb{N}} \mid h \circ f(i_{n})\neq 0\}\) is a finite set. Therefore, \(\{n\in 3{\mathbb{N}} \mid h(i_{n})\neq 0\}\) is a finite set. Let \(n\in{\mathbb{N}}\) such that \(n\neq 3q\) for every \(q\in{\mathbb{Z}}\). Then in \(\prod_{{\mathbb{N}}}^{g}{\mathbb{Z}}\) we have \(i_{n}*i_{n}=i_{n}\). Hence, \(h(i_{n})+h(i_{n})=h(i_{n})\) since \(h\) is a homomorphism. So \(h(i_{n})=0\) since \(A\) is an abelian group. Therefore \(\{n\in{\mathbb{N}} \mid h(i_{n})\neq 0\}\) is a finite set, which implies that \(A\) is a slender generalized group.
\end{proof}

\begin{theorem}
Let \(A\) be an abelian group. If \(A\) is a slender generalized group then it is a slender group.
\end{theorem}

\begin{proof}
Let \(h:\prod_{{\mathbb{N}}}{\mathbb{Z}}\to A\) be a homomorphism of abelian groups. Let \(g:\prod_{{\mathbb{N}}}^{g}{\mathbb{Z}}\to\prod_{{\mathbb{N}}}{\mathbb{Z}}\) be a homomorphism with \(g(x_{1},x_{2},x_{3},x_{4},x_{5},x_{6},\dots)=(x_{3},x_{6},x_{9},\dots)\). Then \(h \circ g:\prod_{{\mathbb{N}}}^{g}{\mathbb{Z}}\to A\) is a homomorphism of generalized groups. Since \(A\) is a slender generalized group, \(\{n\in{\mathbb{N}} \mid h \circ g(i_{n})\neq e(h \circ g(i_{n}))\}\) is a finite set. Since \(A\) is an abelian group, \(e(h \circ g(i_{n}))=0\) for every \(n\in{\mathbb{N}}\). Hence \(\{n\in 3{\mathbb{N}} \mid h \circ g(i_{n})\neq 0\}\) is a finite set. Note that for every \(n\in 3{\mathbb{N}}\), we have \(g(i_{n})=i_{n/3}\). Therefore, \(\{n\in{\mathbb{N}} \mid h(i_{n})\neq 0\}\) is a finite set, which implies that \(A\) is a slender group.
\end{proof}

\section{Further Results on Slender Generalized Groups}

\begin{theorem}
Let \( G \) be a slender generalized group. Then every generalized group homomorphism  
\[
h : \prod_{\mathbb{N}}^g \mathbb{Z} \to G
\]  
has the property that \( h(i_n) = e(h(i_n)) \) for all but finitely many \( n \in \mathbb{N} \).
\end{theorem}

\begin{proof}
This follows directly from Definition 1.6 and the fact that \( G \) is slender.
\end{proof}

\begin{theorem}
Let \( G \) and \( H \) be generalized groups, and let \( f : G \to H \) be a surjective homomorphism. If \( G \) is slender, then \( H \) is also slender.
\end{theorem}

\begin{proof}
Suppose \( h : \prod_{\mathbb{N}}^g \mathbb{Z} \to H \) is a homomorphism. Since \( f \) is surjective, for each \( n \in \mathbb{N} \), there exists \( a_n \in G \) such that \( f(a_n) = h(i_n) \). Define a homomorphism \( \tilde{h} : \prod_{\mathbb{N}}^g \mathbb{Z} \to G \) by \( \tilde{h}(i_n) = a_n \). By the slenderness of \( G \), \( \tilde{h}(i_n) = e(\tilde{h}(i_n)) \) for all but finitely many \( n \). Applying \( f \), we get \( h(i_n) = e(h(i_n)) \) for all but finitely many \( n \). Hence, \( H \) is slender.
\end{proof}

\begin{theorem}
A direct sum of slender generalized groups is slender.
\end{theorem}

\begin{proof}
Let \( \{G_\alpha\}_{\alpha \in I} \) be a family of slender generalized groups, and let \( G = \bigoplus_{\alpha \in I} G_\alpha \). Suppose \( h : \prod_{\mathbb{N}}^g \mathbb{Z} \to G \) is a homomorphism. Then for each \( n \in \mathbb{N} \), \( h(i_n) \) has only finitely many non-identity components. By the slenderness of each \( G_\alpha \), the set of \( n \) for which \( h(i_n) \) has a non-identity component in any \( G_\alpha \) is finite. Hence, \( h(i_n) = e(h(i_n)) \) for all but finitely many \( n \).
\end{proof}

\begin{theorem}
Let \( G \) be a normal slender generalized group. Then every homomorphism  
\[
h : \prod_{\mathbb{N}}^g \mathbb{Z} \to G
\]  
factors through a finite subproduct.
\end{theorem}

\begin{proof}
Since \( G \) is slender, \( h(i_n) = e(h(i_n)) \) for all but finitely many \( n \). Let \( S = \{n \in \mathbb{N} : h(i_n) \neq e(h(i_n))\} \). Then \( h \) is determined by its values on the finite subproduct \( \prod_{n \in S} \mathbb{Z} \).
\end{proof}

\section{Conclusion}
We have extended the classical theory of slender abelian groups to the setting of generalized groups. Our results show that many key properties of slender groups carry over to this more general context, and we have provided a characterization of slenderness in terms of homomorphisms from the generalized product \( \prod_{\mathbb{N}}^g \mathbb{Z} \). Further research may explore the relationship between slenderness and other structural properties of generalized groups, such as topological slenderness in the context of topological generalized groups.

\end{document}